\documentclass[11pt]{article}

\usepackage[nonatbib,preprint]{arxiv}

\usepackage{pdfpages}

\usepackage[utf8]{inputenc} % allow utf-8 input
\usepackage[T1]{fontenc}    % use 8-bit T1 fonts
\usepackage{hyperref}       % hyperlinks
\usepackage{url}            % simple URL typesetting
\usepackage{booktabs}       % professional-quality tables
\usepackage{amsfonts}       % blackboard math symbols
\usepackage{nicefrac}       % compact symbols for 1/2, etc.
\usepackage{microtype}      % microtypography
\usepackage{lipsum}
\usepackage{graphicx}
\graphicspath{ {./images/} }

\usepackage{times}

\usepackage{amsfonts,amssymb,amsmath,amsthm}
\usepackage{wrapfig}
\usepackage{latexsym}
\usepackage{psfrag}
\usepackage{multirow}
\usepackage{mathrsfs}
\usepackage{color}
\usepackage[cmtip,arrow]{xy}
\usepackage{pb-diagram,pb-xy}
\usepackage{graphicx}

\usepackage{tikz}
\usetikzlibrary{matrix,arrows,positioning,calc,shapes}
\usepackage{graphicx}
\usepackage{subfigure}
\newcommand{\sA}{{\mathsf A}}

\newcommand{\sJ}{{\mathsf J}}

\newcommand{\sM}{{\mathsf M}}

\newcommand{\R}{{\mathbb{R}}}

\newcommand{\cF}{{\mathcal F}}
\newcommand{\cX}{{\mathcal X}}

\newcommand{\Int}{\mathop{\mathrm{int}}\nolimits}

\newcommand{\mvmap}{\rightrightarrows}

\newcommand{\sABlock}{{\mathsf{ABlock}}}
\newcommand{\pred}[1]{\overleftarrow #1}
\newcommand{\supp}[1]{\left| {#1}\right|}
\newcommand{\setof}[1]{\left\{ {#1}\right\}}

\title{Extracting Global Dynamics of Loss Landscape in Deep Learning Models}

\author{
 Mohammed Eslami \\
  Netrias, LLC\\
  \texttt{meslami@netrias.com} \\
     \And
Hamed Eramian \\
  Netrias, LLC\\
  \texttt{eramian@netrias.com} \\
   \And
Marcio Gameiro \\
  Rutgers University\\
  \texttt{gameiro@math.rutgers.edu} \\
     \And
Konstantin Mischaikow \\
  Rutgers University\\
  \texttt{mischaik@math.rutgers.edu} \\
     \And
William Kalies \\
  Florida Atlantic University\\
  \texttt{wkalies@fau.edu} \\
}

\begin{document}
\maketitle
\begin{abstract}
Deep learning models evolve through training to learn the manifold in which the data exists to satisfy an objective. It is well known that evolution leads to different final states which produce inconsistent predictions of the same test data points. This calls for techniques to be able to empirically quantify the difference in the trajectories and highlight problematic regions. While much focus is placed on discovering what models learn, the question of how a model learns is less studied beyond theoretical landscape characterizations and local geometric approximations near optimal conditions. Here, we present a toolkit for the Dynamical Organization Of Deep Learning Loss Landscapes, or DOODL3. DOODL3 formulates the training of neural networks as a dynamical system, analyzes the learning process, and presents an interpretable global view of trajectories in the loss landscape. Our approach uses the coarseness of topology to capture the granularity of geometry to mitigate against states of instability or elongated training. Overall, our analysis presents an empirical framework to extract the global dynamics of a model and to use that information to guide the training of neural networks. 
\end{abstract}

% keywords can be removed

\section{Introduction}
Deep learning models evolve through training to learn the manifold in which the data exists to satisfy an objective. The evolution of these transformations often does not follow a consistent path for every initial condition and so predictions of the same test points can vary significantly with no change in model architecture or hyperparameters. This is due to a well-known phenomenon that different final states are reached given a different set of initial weights when a model is trained and tested with the same data and hyperparameters. While there is a significant amount of work in theoretically characterizing global properties of the loss landscape of a neural network, an empirical analysis of the global dynamics that emerge as a result of training is less studied and understood. We define global dynamics as the trajectory of the family of models given a set of parameters. This is in contrast to the often presented local analysis that focuses on general properties or loss evaluation near optima. 

Identification and characterization of regions of dynamic sensitivity is seldom studied in the deep learning domain, but well studied in nonlinear dynamics mathematics. This is due to the fact that global analysis of dynamical systems with a large parameter space ($>10$) is often computationally intractable. Deep learning models can contain millions to billions of parameters and thus are not a good fit for these methods. However, discovery and mitigation of regions of sensitivity is a requirement to build models that generalize and need an analysis of the training dynamics. Formulation of the training process of a deep learning model as a dynamical system and modeling the system's dynamical properties can enable such analyses. Here, we present a toolkit for the Dynamical Organization Of Deep Learning Loss Landscape, or DOODL3. DOODL3 uses techniques from dynamical systems theory to extract the global dynamics of training a model and empirically identifies regions of sensitivity, elongated training, or stable attractors. Our approach exploits the malleability of topology with the goal of providing concrete guidelines towards the construction of geometric models and proofs for specific systems. It uses the coarseness of topology to enable the us to capture the granularity of geometry that can be used to mitigate against states of instability or elongated training. To the best of our knowledge, this is among the first empirical analyses to globally characterize the dynamics of training neural networks. Thus, we present a The theoretical framework that justifies the claim that topological approximations of training dynamics are sufficient to accurately discover the loss landscapes. Then, we present the DOODL3 toolkit and experimental results that test the relevance of this approach for the deep learning community.

\section{Related Work}

The characterization of a deep learning model's landscape is not a new concept and has been well studied for convergence, stability, and robustness properties \cite{Landscape_nips01,global_convergence_optimal_transport, goodfellow2015qualitatively,chizat2018global, viz_loss_landscape}. Many of the recent advances, however, rely on the theoretical foundations to discover particular properties of the landscape or empirically analyze small regions of the landscape to determine its impact on stability of the model. Given a characterization, \cite{Garipov_2018} presents a method that interpolates between multiple modes in a 2D space to create a low loss tunnel. This work is extended empirically by \cite{Fort_2019}, whose core goals align with those presented in this paper, to model the landscape as a set of interconnected wedges that form an interconnected structure. They show subspaces can connect a set of solutions but do not present the most likely solution given a set of initial conditions. Other related work chooses to focus on the characterization of boundaries between classes \cite{boundary_char} and uses those properties to generate additional data for further characterization. Finally, and most recently, \cite{Deep_Ensembles} empirically identifies diverse modes in function space of the system but does not provide a global characterization of the landscape. Their diversity-accuracy metric is an invaluable metric to quantify a system's predictive shift and is a metric we adapt for our work as well. They present strategies to construct a subspace of the trajectories and analyze the similarity between the predictions. In all of these works, emphasis is placed more on the structure of the landscape over the dynamics introduced by the training process. Most relevant to our work is the work of \cite{viz_loss_landscape}, where they extract the chaotic nature of the loss landscape in 1D, 2D, and 3D surfaces as well as conduct an analysis of specific trajectories through those surfaces. 

The impact of a global understanding of the landscape can also be understood from existing work. 

\begin{itemize}
    \item Reduced Generalization Error: It has been recently shown \cite{GoogleUnderspec} that underspecification of models with a large number of parameters compared to the number of data points can lead to inconsistent predictions. Our analysis of global dynamics can quantify the sensitivity of a model to initial conditions, thereby quantifying the extent of error. 
    \item Control for Faster Learning: Global views provide an escape plan for model parameters if regions of instability are unavoidable. It can further offer a potential to reach specific, desired regions of stability \cite{seidman2020robust}.
    \item Stable Regions of Transfer: Guarantees over families of models will enable effective transfer to new architectures and tasks \cite{wang2018theoretical}.
\end{itemize}

Our goal is to build upon the foundations set in the characterization of the loss landscape with a suite of tools that can empirically characterize training dynamics to accurately discover the loss landscape.

\section{Novel Application of Conley Theory for Characterizing Loss Landscapes of Neural Networks}
\label{sec:background}

DOODL3 focuses on the characterization of dynamics within the loss landscape of a training map for classification functions that are obtained by using gradient descent like methods.
Conceptual images of these loss landscapes typically consist of a one dimensional curve with multiple critical points or an undulating two dimensional surface in three dimensions, and the local minima, or ideally global minimum, represent a final state  of the training process. 
Of course, as indicated in the introduction, the space of parameters $X$ in current deep learning models can be very high dimensional.
It is important to realize that in general one cannot hope to generate or store a geometrically precise representation of such high dimensional surfaces. DOODL3 utilizes a framework to compute an alternative characterization.

The need to characterize the structure of high dimensional spaces is not a novel problem.
A practical solution is to replace the goal of a geometric description with an understanding of the global topology.
This is often achieved using algebraic topology, e.g.\ homology.
The reader may be familiar with the setting of smooth manifolds where \emph{Morse theory} provides the appropriate tools \cite{Morse}.
A compact Riemannian manifold $M$ and a function $f\colon M\to \R$ satisfying certain generic conditions gives rise to a gradient flow where orbits are either isolated equilibria (critical points of $f$) or heteroclinic trajectories limiting in forward and backward time to the equilibria. 
The homology of $M$, denoted by $H_*(M)$, can then be recovered from the \emph{Morse indices} of the critical points and heteroclinic trajectories between equilibria whose Morse indices differ by one. 
A useful summary is that Morse theory is based on two sets of information.
First, there is a combinatorial/order theoretic description of the dynamics, i.e.\ an acyclic directed graph
(equivalently a partially ordered set or poset)
where the vertices are the equilibria and the edges indicate the existence of a heteroclinic orbit from an equilibrium of Morse index $n$ to an equilibrium of Morse index $n-1$.
Second there is algebraic information in the form of a \emph{chain complex}, i.e.\ matrices whose entries indicate the number (with appropriate signs) of heteroclinic orbits between the equilibria of appropriate index.

The mathematical requirements of Morse homology are too strong to be directly applied to training in the context of machine learning.
In practice, the process of training involves a map not a flow, the set of orbits can be much more complicated than isolated fixed points and heteroclinic connections (in fact because of overparameterization the fixed points need not be isolated), and it cannot be assumed that the particular training map employed satisfies the above mentioned generic conditions.
Thus, we turn to \emph{Conley theory}, which is a topological version of Morse theory.
We explain this theory by first focusing on the order theory and second briefly mentioning the algebra.  

\subsection{Conley-Morse graphs}\label{subsec:Conley-order}

 We start our description of Conley theory with the concept of an \emph{attracting block}, defined as a compact set $A\subset X$ such that under the dynamics $A$ is immediately moved into $\Int(A)$, the interior of $A$.
The set of all attracting blocks, $\sABlock,$ forms a lattice \cite{kalies:mischaikow:vandervorst:14}, i.e.\ a poset ordered by inclusion with the property that the intersection or union of attracting blocks is an attracting block. 
In general $\sABlock$ has uncountably many elements, however our use of Conley theory is based on identifying a finite sublattice $\sA\subset \sABlock$.

For computational purposes, we first discretize the phase space $X$ to obtain a finite decomposition.
The theoretical constraints on the decomposition are quite weak, but from a practical perspective it is reasonable to assume that $X$ is decomposed into polytopes that are indexed by a set $\cX$.
In particular, given $\xi\in\cX$ we let $\supp{\xi}\subset X$ denote the associated polytope.
We also assume that the dynamics, e.g.\ training map, is given by a continuous function $f\colon X\to X$.
Ignoring for the moment the computational cost, assume that for each polytope $\xi$ we can determine $\cF_0(\xi):= \setof{\xi'\in\cX\mid f(\supp{\xi}|)\cap \supp{\xi'}\neq \emptyset}$.
Observe that we have defined a multivalued map $\cF_0\colon \cX \mvmap \cX$, that can either be interpreted as a directed graph with vertices $\cX$ and edges $\xi\to\xi'$ if $\xi'\in \cF_0(\xi)$, or as a map from polytopes to sets of polytopes.
A multivalued map $\cF\colon\cX\mvmap \cX$ is called an \emph{outer approximation} of $f$ if $\cF_0(\xi)\subseteq \cF(\xi)$ for all $\xi\in \cX$. The computations on neural nets in Section~\ref{sec:results} is performed using multiscale cubical grids as shown in Figure~\ref{fig:multi_scale}.

Combinatorial multivalued maps are used to characterize global dynamics at a fixed resolution, 
and our primary interest is to identify the recurrent versus nonrecurrent dynamics.  
Viewed from the perspective of a directed graph, this is a question of identifying the strongly connected components of $\cF$ for which there are extremely efficient algorithms, such as Tarjan's algorithm.
A \emph{recurrent component} is a strongly connected component that contains at least one edge.
The condensation graph of $\cF$ is obtained by identifying each strongly connected component to a single vertex. 
Note that this is a directed acyclic graph and hence a poset $\mathsf{SC(\cF)}$.
The \emph{Morse graph} of $\cF$, denoted by $\sM(\cF)$, is the 
subposet consisting of recurrent components of $\cF$. For an outer approximation, the recurrent dynamics of $f$ is contained in the regions corresponding to elements of  $\sM(\cF)$ so that outside of the recurrent components of $\cF$, the dynamics of $f$ is nonrecurrent.

The minimal elements in $\sM(\cF)$ correspond to attracting blocks in $X$ in which one could expect local minima of the loss function to reside. Non-minimal recurrent components may contain unstable critical points of the loss function. The nonrecurrent region contains basins of attraction for the minimal recurrent components as well as boundaries between these basins or \emph{separatrices}. To identify these basins and separatrices, we look for an \emph{order retraction of  $\mathsf{SC(\cF)}$ onto $\sM(\cF),$ } i.e.\ an order-preserving map $\sigma\colon\mathsf{SC}(\cF)\to\sM(\cF)$ such that $\sigma\circ i=\text{id}_{\sM(\cF)}.$ An order retraction need not exist, but when present, it characterizes a consistent organization of the global dynamics into minimal recurrent sets, their basins of attraction, and their separatrices. An algorithm for computing an order retraction, or that one does not exist, is given in \cite{KKV}. If no order retraction exists, this is an indication that the dynamics is not well resolved by $\cF$.

Thus, at least conceptually we can achieve the goal of providing a combinatorial description of the loss landscape: the elements of $\sM(\cF)$ capture potential recurrent dynamics and the partial order on $\sM(\cF)$ identifies the nonrecurrent dynamics.
It is important to observe that the computations are being done using the outer approximation $\cF$, but the mathematical results apply to the object of interest, the dynamics of the continuous map $f$.

Focusing in on the problem of machine learning, given that $X$ is high dimensional and there is limited training data, it is unreasonable to expect that $f$ is known for all of $X$.
However, as is emphasized above we do not need precise knowledge of $f$ for the mathematical machinery to be employed.
In particular, let us assume the following information $\setof{(x_i,f(x_i)) \mid x_i \in X, i=1,\ldots, I }$.
We can use this to produce a surrogate model $F\colon X\to X$.
Given $F$, fix a variance $\mu$ and determine $\cF$ an outer approximation to this variance, i.e.\ the image of $\supp{\cF(\xi)}$ captures the image of $F$ within the variance over $\supp{\xi}$.
Applying the procedure described above will produce a Morse graph $\sM(\cF)$ and the variance provides a lower bound on the probability that $\sM(\cF)$ is a correct poset representation for the global dynamics of $f$.

Returning to lattices of attracting blocks, observe that associated to each $M\in \sM(\cF)$ is a set of elements of $\cX$ and hence polytopes in $X$.
Let $A(M)$ denote the set of polytopes that are reachable from $M$ via an outer approximation $\cF$.
Then $\setof{A(M)\mid M \in\sM(\cF)}$ generates a finite sublattice  $\sA(\sM(\cF))\subset\sABlock(f)$, see
\cite{kalies:mischaikow:vandervorst:15,KKV,kalies:mischaikow:vandervorst:21}.
Therefore given $\cF\colon\cX\mvmap \cX$ we have an extremely efficient construction of the building blocks of Conley theory.
Finally, in \cite{kalies:mischaikow:vandervorst:15,kalies:mischaikow:vandervorst:21} it is shown that given a Conley-Morse graph for the dynamics of $f$ there is no theoretical obstruction to realizing it via $\sM(\cF)$ for an outer approximation $\cF$ other than computational effort.

\subsection{Conley index}\label{subsec:Conley-index}
Similar to Morse theory, Conley theory also has an algebraic component, the Conley index, which we now briefly describe. 
An element of a bounded, distributive lattice $\sA$ is  \emph{join irreducible} if it has a unique immediate predecessor $\pred{A}$ with respect to the partial order of inclusion.
We denote the set of join irreducible elements of $\sA$ by $\sJ(\sA)$.
In the setting of a flow if $A\in \sA\subset\sABlock$ is join irreducible, its \emph{Conley index} is $H_*(A,\pred{A})$.
The associated chain complex is given by linear maps between the Conley indices of elements of $\sJ(\sA)$ \cite{harker:mischaikow:spendlove}.
It is important to note that at this step $\sJ(\sA)$ is being used as an indexing set, and thus, it is  the poset structure of $\sJ(\sA)$ that plays a role. When an order retraction exists in the computations described in the previous subsection, Birkhoff's theorem implies that $\sJ(\sA(\sM(\cF)))$ is isomorphic to $\sM(\cF)$ so that the poset structure is identified.

The take away message is that in Conley theory  $\sJ(\sA)$ plays the role of equilibria in Morse theory and knowledge of partial order on $\sJ(\sA)$ corresponds to knowledge about heteroclinic orbits. 

The focus of this paper is on the identification of $\sJ(\sA)$ in the context of training a neural network.
We note, as this is relevant for future work, in both settings linear maps provide information about the existence of trajectories between invariant sets. The additional feature of Conley theory is that the Conley index provides information about structure of the dynamics in $A\setminus \pred{A}$, which, as was emphasized above is expected to be more complicated than an isolated equilibrium. 

With respect to the computational methods that are the focus of this paper, 
a powerful feature of Conley theory is that we obtain an accurate characterization of the global dynamics by correctly identifying the poset structure of  $\sJ(\sA)$ and $H_*(A,\pred{A})$ for each $A\in \sA$.
In particular, since homology is a homotopy invariant one can expect to correctly compute the Conley index even if the sets in $\sA$ are not attracting blocks for the dynamics. This is the case in the setting of this paper where machine learning models are studied using a finite set of training data and surrogate models as described in the previous subsection.

\section{DOODL3 Toolkit}

The theoretical framework discussed in Section~\ref{sec:background} justifies the claim that topological approximations of training dynamics are sufficient to accurately discover the loss landscapes. However, the theory does not guarantee that the application of this method will necessarily lead to interesting discoveries, e.g.\ it is possible that output is a single attracting block and all points in parameter space map to that block. We developed DOODL3 to test its relevance for the deep learning community.

DOODL3 takes a dataset, a neural network architecture, an optimizer, and a loss function as input and will analyze the neural networks training dynamics. In its current implementation, DOODL3 varies initial weights using a uniform distribution in parameter space. The neural network is then trained with a seeded train/test split of the dataset for $N$ epochs. This process is repeated $K$ times and checkpoints output the initial and final weights of each cycle. Note that to discover the global training dynamics, DOODL3 does not need the intermediate steps of the trajectory. DOODL3 takes a collection of initial set of weights, $\{ X^0_k \mid k = 1, \ldots, K \}$ sampled to cover a large region of parameter space, and a collection of the final set of weights, $\{ X^N_k \mid k = 1, \ldots, K \}$, and trains a surrogate map, $F$, a Gaussian process, to approximate the training dynamics. 

\begin{figure}[hbt]
\includegraphics[width=12cm, height=2.5cm]{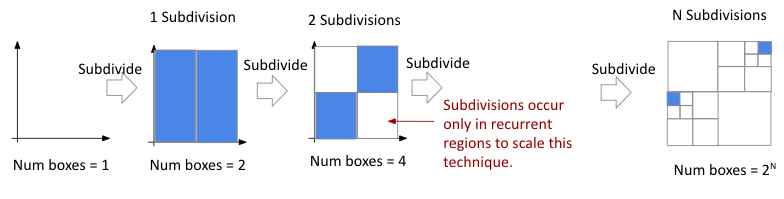}
\centering
\caption{Multiscale Adaptive Mesh-based Decomposition of Phase Space}
\label{fig:multi_scale}
\end{figure}

To construct a multivalued map $\cF$ we start with an adaptive mesh based decomposition of the parameter space as described in \cite{bush:gameiro:harker:kokubu:mischaikow:obayashi:pilarczyk} and shown in Figure \ref{fig:multi_scale}, and then use $F$ to sample across the grid. Finally, as indicated in Section \ref{subsec:Conley-order} we compute $\sM(\cF)$ and  an order retraction. The Morse graph and their underlying parameter regions are further accessible through APIs for further analysis. 

\section{Experiments and Results}\label{sec:results}

\begin{figure}[hbt]
\centering
\includegraphics[width=1\linewidth]{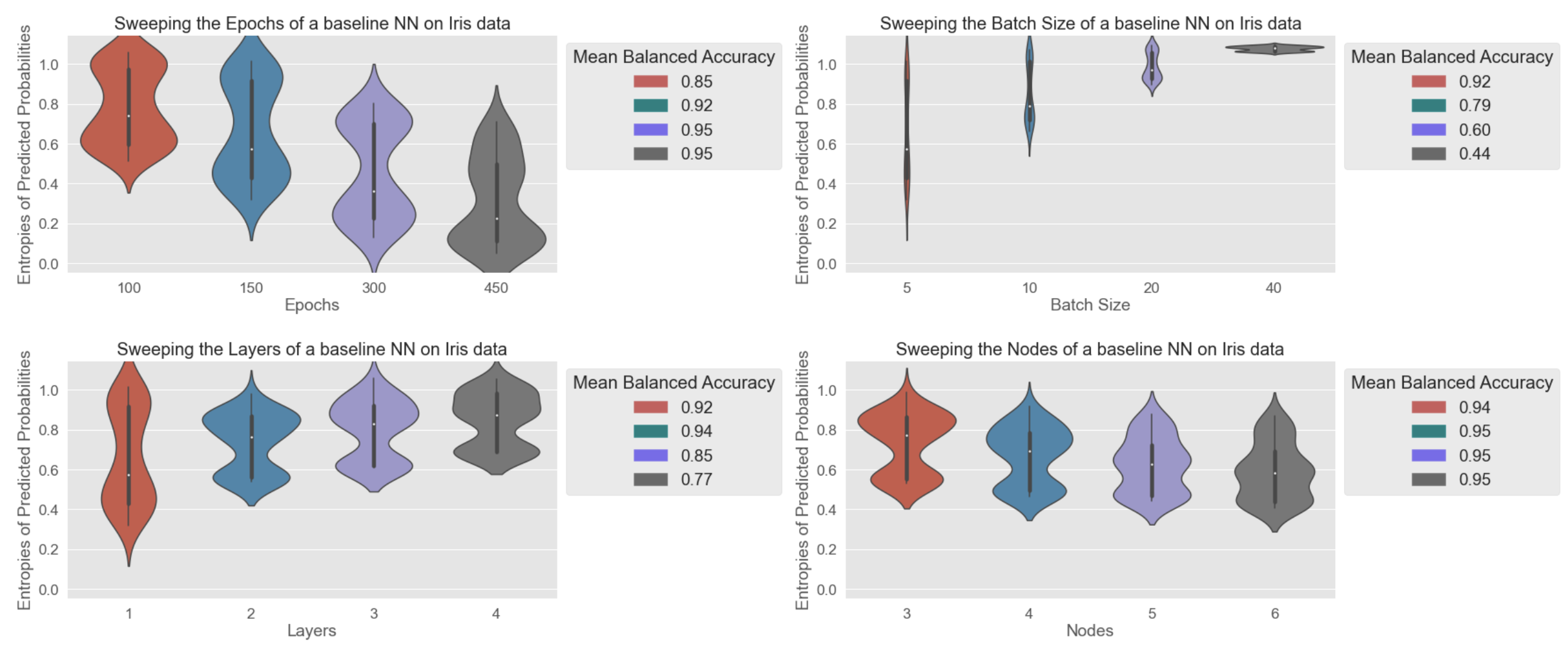}
\caption{ Distribution of entropy across every test point in iris dataset. A single layer, single node neural network trained for 150 epochs was used as the baseline. X-axis indicate additional hyperparameters and their ranges that were tested.}\label{fig:entropy_iris}
\end{figure}

Following \cite{Fort_2019}, we tested DOODL3 on two toy examples, a fully connected neural network (FCNN) and a convolutional neural network (CNN). The FCNN used the Iris dataset, trained a baseline model, and then swept the number of layers, number of nodes, epochs, and batch size to generate training dynamics information. Each sweep had 100 training cycles and we evaluated the entropy of each test point across all cycles (Fig.~\ref{fig:entropy_iris}). While the average balanced accuracy is high, a wide distribution with respect to the entropy for each condition exists. As accuracy increases, the average entropy decreases but the distribution remains wide. Thus, this model is a good candidate for analysis of dynamics. 

\begin{figure}[hbt]
\centering
\includegraphics[width=1\linewidth]{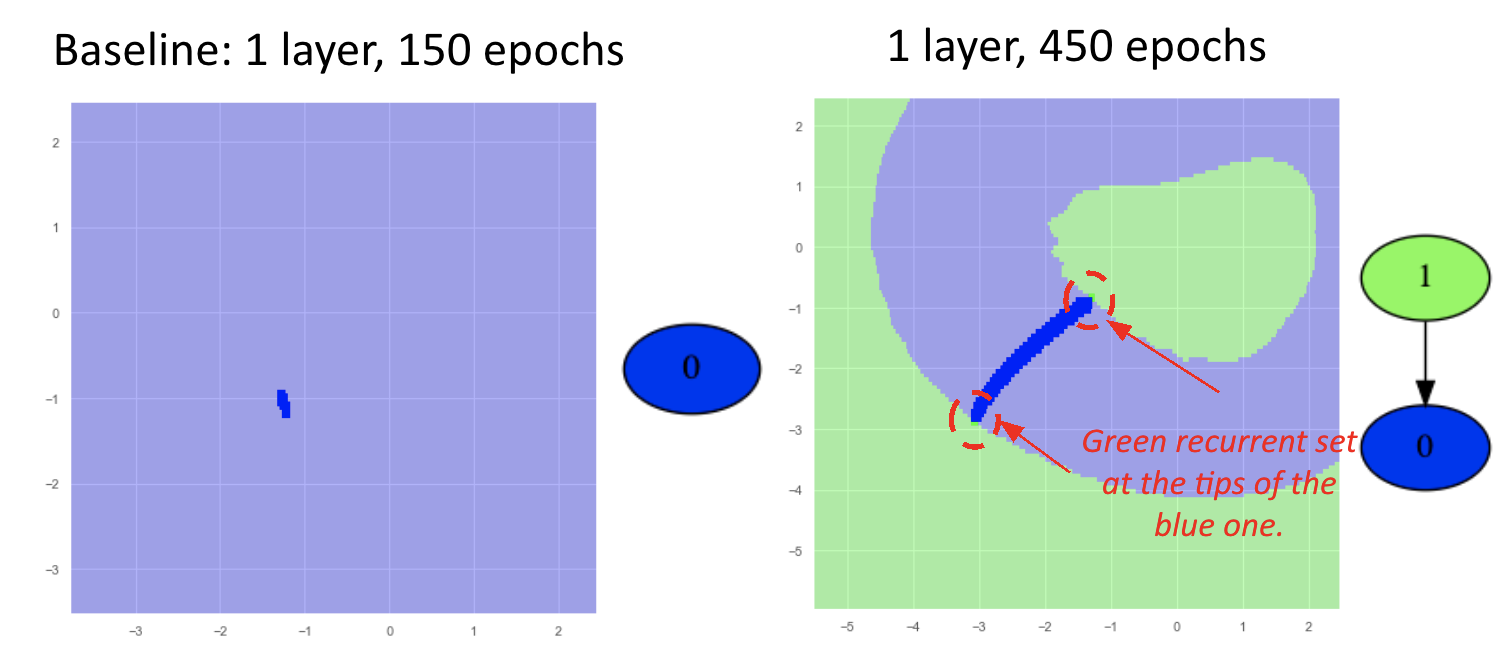}
\caption{Global dynamics of model training during for the iris dataset.}
\label{fig:dynamics_iris}
\end{figure}

An analysis of the weights showed that a majority of the weights in the network did not change significantly. Thus, we selected two weights and a set of conditions whose entropy distributions were wide but accuracies were low to compare to a baseline model (Fig.~\ref{fig:dynamics_iris}).

With an increased number of epochs the Morse graph changed from a single node to two ordered nodes (Fig.~\ref{fig:dynamics_iris}(left) vs. (right)).

In the baseline case, the detected dynamics presents as a small global attractor whose basin of attraction contains the entire rectangular region. With increased epochs, the recurrent region has elongated and separated into to two components. The Morse graph and order retraction show that these two components can be grouped together again as a global attractor for the entire region. However, trajectories that begin in the green region may exhibit different transient behavior, perhaps requiring additional iterations. The Morse graph and the topology of the corresponding regions provides information about the loss landscape and training dynamics. 

\begin{figure}[hbt]
\centering
\includegraphics[width=1\linewidth]{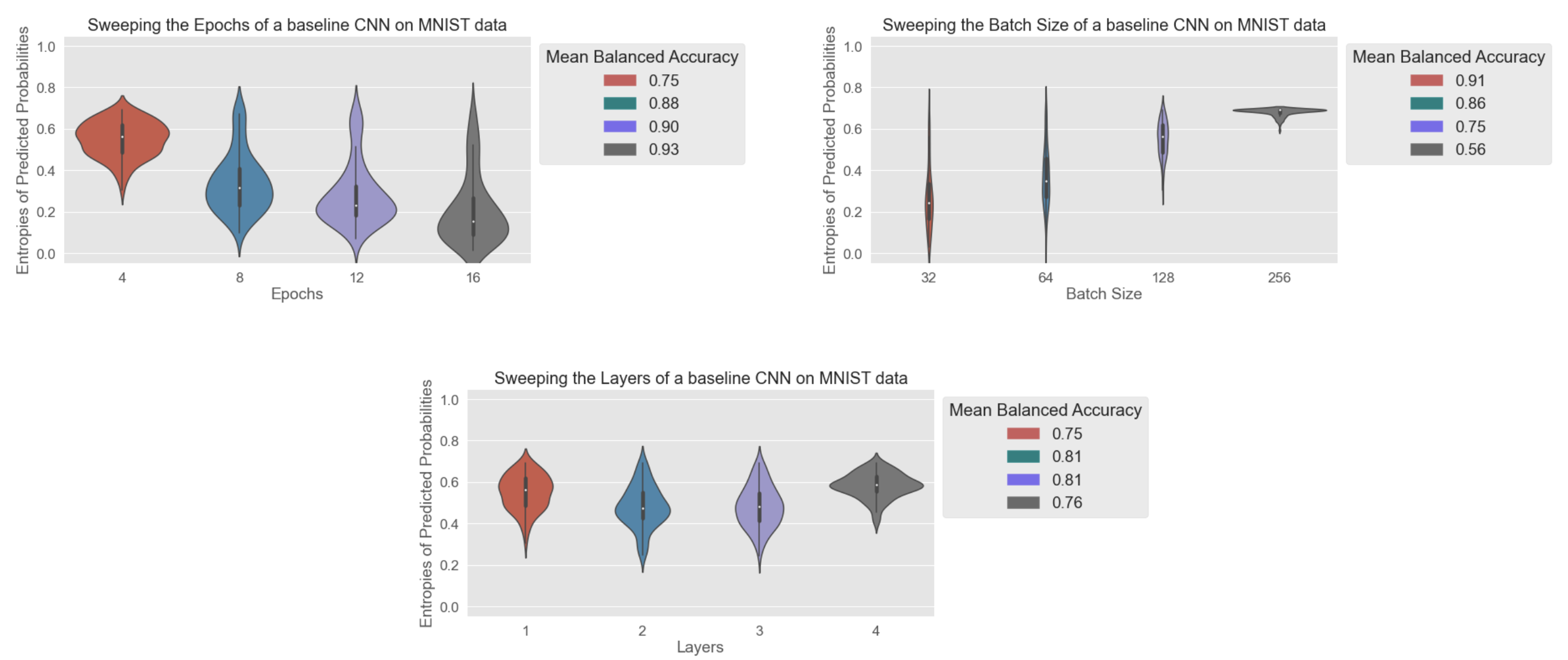}
\caption{Distribution of entropy for the MNIST dataset. A single convolutional layer connected to a soft-max layer trained for 4 epochs was used as the baseline model. X-axis indicate additional hyperparameters and their ranges that were tested.}
\label{fig:entropy_mnist}
\end{figure}

\begin{figure}[hbt]
\centering
\includegraphics[width=12cm, height=4cm]{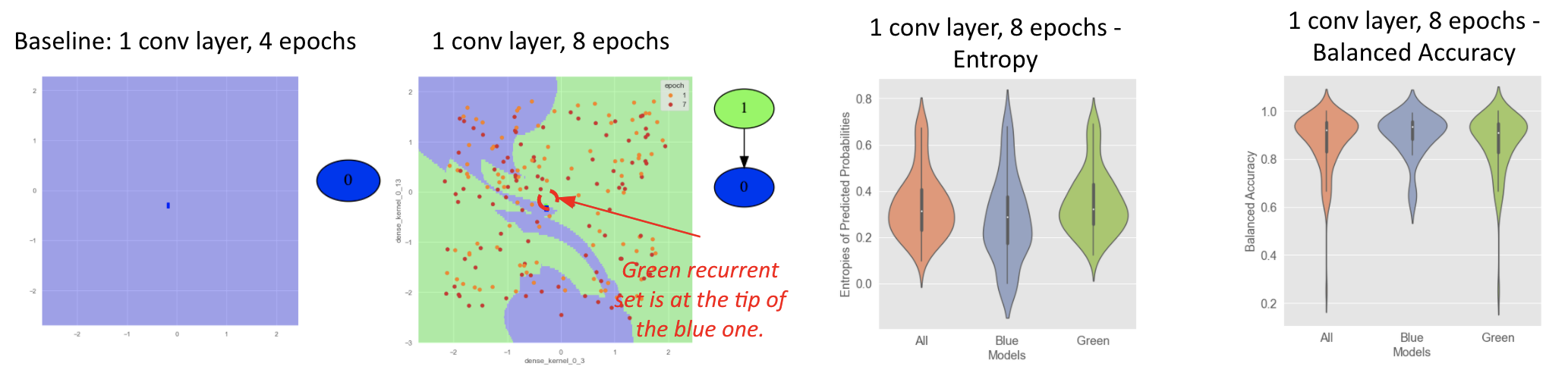}
\caption{Global dynamics of model training during for the MNIST dataset}
\label{fig:dynamics_mnist}
\end{figure}

The CNN was  trained with the MNIST dataset. Again, the entropy of the test points and average accuracy showed that this model was a good candidate for dynamics analysis (Fig.~\ref{fig:entropy_mnist}). As was the case with the FCNN, the global dynamics showed an increase in the number of recurrent sets as the number of epochs increased. We hypothesize that this is due to the fact that the additional epochs provide additional information to the multivalued map to refine the dynamics approximations (Figure \ref{fig:dynamics_mnist}). Furthermore, more interestingly, an analysis of predictions made by the model showed that points whose initial conditions that were mapped to the blue recurrent set were able to achieve lower entropy values as well as a more tightly distributed balanced accuracy. The subtle difference, we believe, is due to the proximity of the green recurrent set to the blue, but, nonetheless a difference is highlighted.

\section{Conclusions}

DOODL3 exploits the malleability of topology with the goal of providing concrete guidelines towards training neural networks. The tools in DOODL3 are based on the minimal perspective of characterizing the global dynamics of training a neural network with recurrent and non-recurrent sets. It uses the coarseness of topology to capture the granularity of geometry (regions of stability) to discover how a model learns. Finally, the Morse graph decomposition provides an interpretable view of the high-dimensional global dynamics for further analyses. 
While significant effort has been placed on scaling this analysis, applying this analysis on neural networks with millions of parameters will require more work to improve the efficiency. First, moving from a mesh-based decomposition of parameter space to mesh-less one can significantly improve the extraction of dynamics. Second, techniques to  adaptively identify important regions coupled with coarse dynamic representations can inform the sampling of the space for reconstruction of the model’s dynamics.

\section{Broader Impact}

Artificial intelligence, and deep learning in particular, is increasingly becoming a critical component in our national infrastructure--- from scientific discovery to automation and manufacturing, to fraud detection and cybersecurity. It is thus imperative that techniques used to train these models be consistent in their predictions. DOODL3 aids in identifying and characterizing regions of sensitivity, as well as the ability to control the sensitivity to build generalizable, safe, and consistent AI models. DOODL3 can be applied to any neural network architecture with any choice of underlying optimizer.

\section{Acknowledgements and Disclosure of Funding}
This material is based upon work supported by the United States Air Force, Air Force Office of Scientific Research, Arlington VA and the Defense Advanced Research Project Agency (DARPA) under Contract No. FA9550-20-C-0001. Any opinions, findings and conclusions or recommendations expressed in this material are those of the author(s) and do not necessarily reflect the views of the United States Air Force, Air Force Office of Scientific Research, Arlington VA and DARPA.

\bibliographystyle{unsrt}  
\bibliography{references}  



\section*{Supplementary material}

# S1. Toolkit Availability and Reproducibility

The code, data, and instructions to generate the plots in the paper are located here: https://github.com/netrias/DOODL3 If one is interested in only constructing a Conley-Morse Graph database (CMGDB), the tools there are located here: https://github.com/marciogameiro/CMGDB

Given that DOODL3 builds its corpus of data using a deep learning model using the Adam optimizer, a stochastic dynamical system at its core, we should highlight that the results may not be identical to the ones highlighted in the paper. In an effort to allow readers to use from the exact dataset used in the paper, we have placed the data used for this paper in a data folder in in the DOODL3-master directory.

Furthermore, given the creation of Conley-Morse graphs is a numerical process that depends not only on the surrogate map (see S2 for details) but also on extracting connected components using Tarjan's algorithm, the outputs can lead to slightly different structures but the global properties are preserved. Figure S.1 is an example of different outputs for the iris dataset that were received using the same settings and libraries:

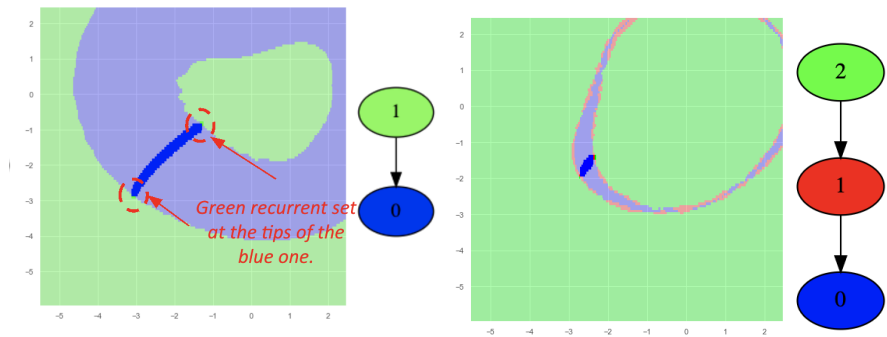

Figure S.1: Morse graphs and 2D Order retractions of the iris dataset using the same settings and hyperparameters for calculations.

While the two are different, they both have the same coarse topologies and thus will not affect downstream proofs that will depend on them.

Another parameter that can be changed is the number of subdivisions, which is the number of times the grid is adaptively refined as discussed in Section 3 of the paper. This will allow one to increase the resolution of the extracted topology to further dissect recurrent sets that are identified. The changes in these settings, as expected, are more apparent on the output.

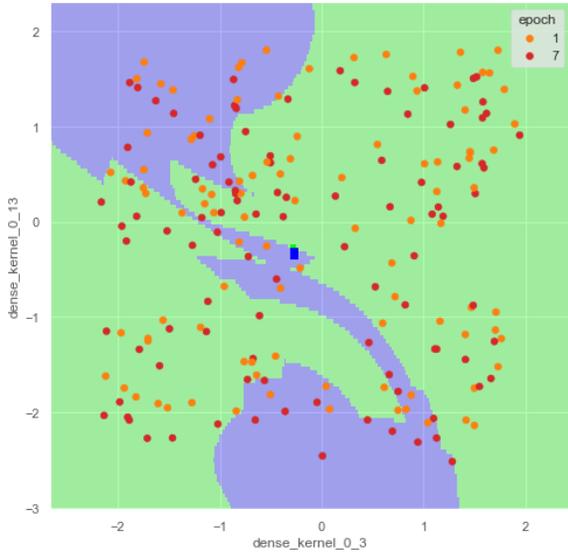

(a)

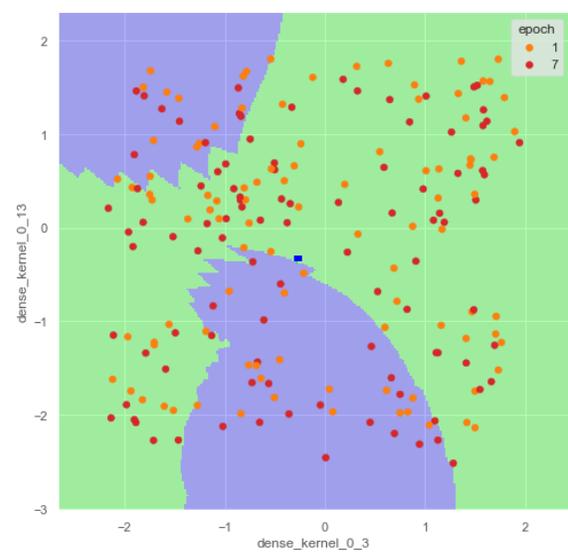

(b)

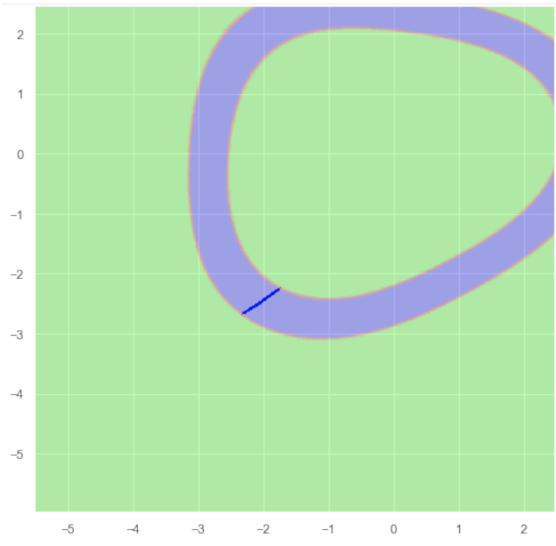

(c)

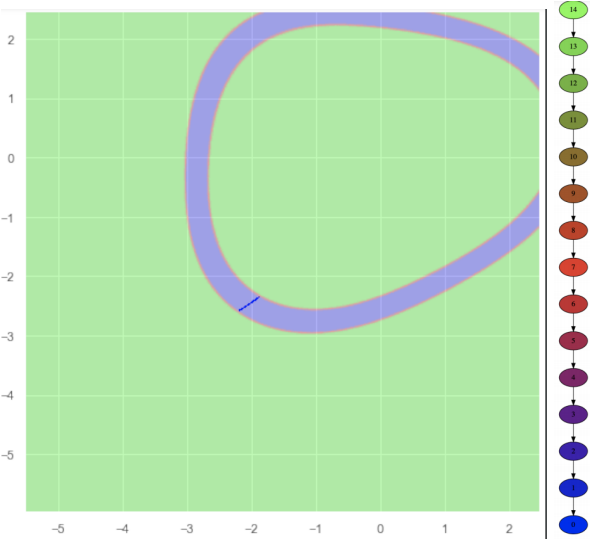

(d)

Figure S.2. Order retraction and associated Morse graphs (if they are different from what was shown in the manuscript) for the (a) convolutional network at 15 subdivisions, (b) convolutional network at 17 subdivisions, (c) fully connected network at 22 subdivisions and (d) fully connected network at 24 subdivisions.

The convolutional network refines the 2D projection that was extracted with 15 subdivisions but its Morse graph is preserved to what was shown in the primary manuscript. The iris dataset, however, shows a continued refinement of the separaptrix between the green and purple regions. The overall topology remains the same but we need to conduct additional tests to see if the additional Morse nodes identified are indeed disjoint recurrent sets or a case of overfitting.

2

# S2. Computing Entropy and Accuracy Metrics Across The Model Runs

A model's consistency can be measured in terms of the predictions it makes on the same test data points across independently perturbed training instances of the model. While prior work uses an aggregate measure of the percent of test point labels that have changed, we choose a measure that reflects the consistency of each test point. We believe a metric per datapoint can disambiguate aggregate statistical measures that could show points near classification boundaries as being sensitive. An analysis of the final states and predictions of the fully connected network shows that models with random initial weights that are near each other can make vastly different predictions and those points that are not consistently predicted are often not boundary points (Figure S.3)

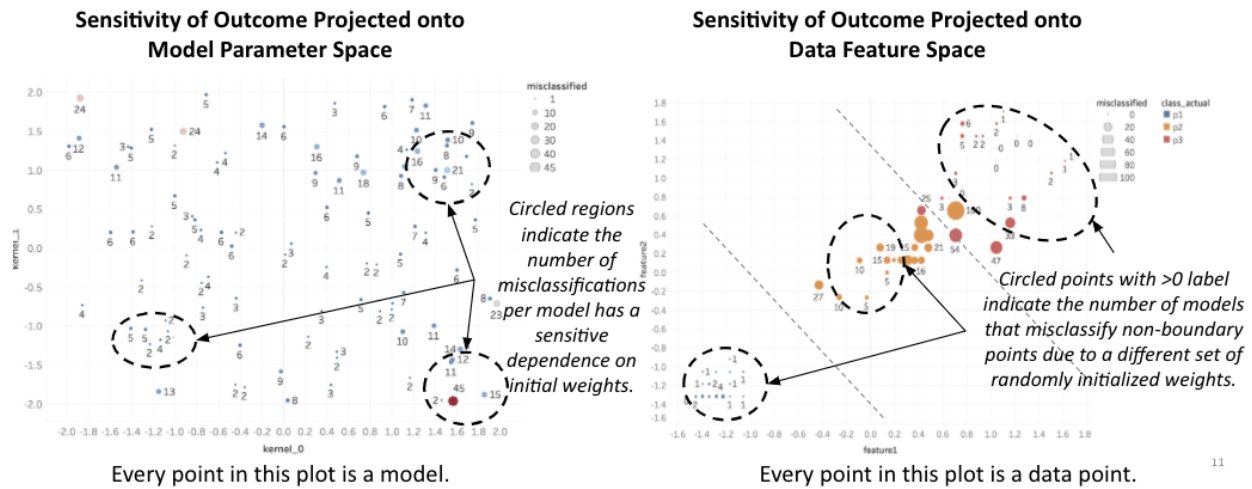

Figure S.3. Order retraction and associated Morse graphs (if they are different from what was shown in the manuscript) for the (a) convolutional network at 15 subdivisions, (b) convolutional network at 17 subdivisions, (c) fully connected network at 22 subdivisions and (d) fully connected network at 24 subdivisions.

Given the entropy distributions in the primary text (Figure 2 and Figure 4), we see that the model will have points that are difficult to predict, namely ones with high entropy. The natural next question then is whether these points have an impact on the performance of the models? Namely, can the model overcome data points with high-entropy to achieve high accuracies, even if it is not confident. If so, then only using accuracy as a metric may yield blind spots where the model is not confident, but correct, and a small perturbation could change its predictive performance. Ideally, we want models whose performance improvement also increases its confidence. In Figure S.4 we see distributions of balanced accuracies from Iris models trained with different hyperparameters. The 100 models were trained and tested for each set of hyperparameters. We notice that on average, the distributions of balanced accuracy converge to much tighter distributions than distributions of entropy seen in Figure 2 of the main text. Thus,

while the distributions of performance are tight, there is a wide range of confidence in those predictions.

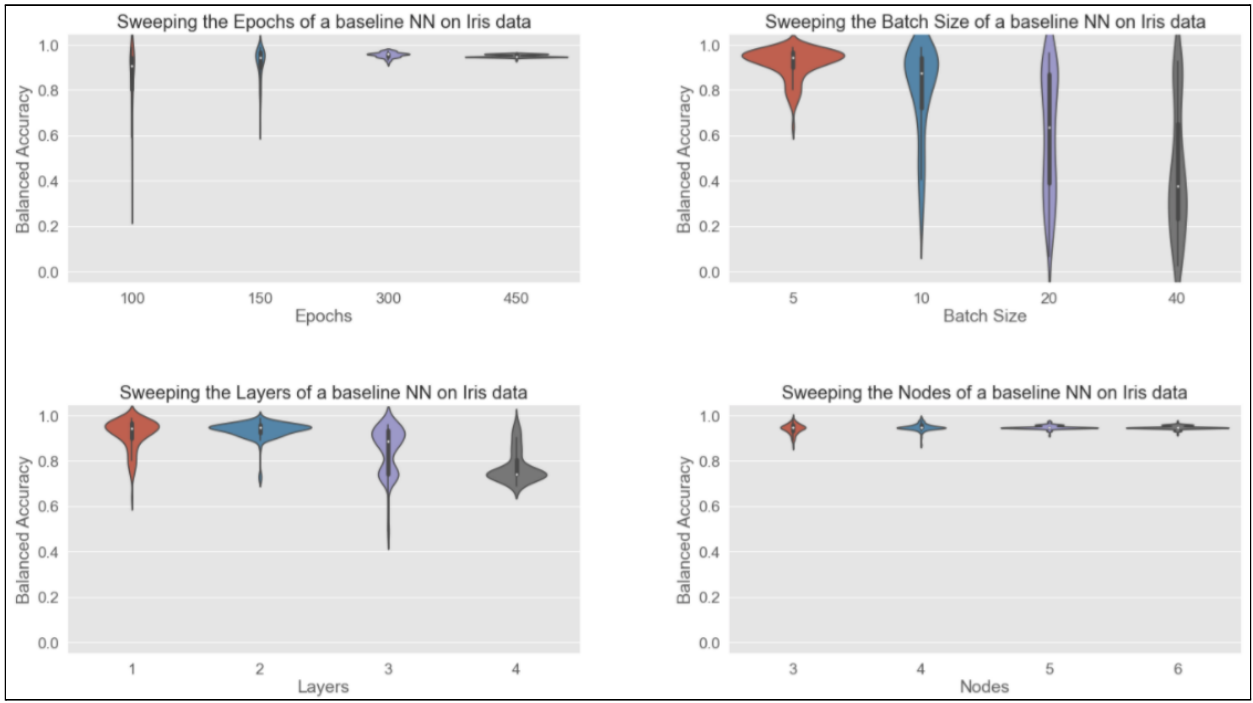

Figure S.4. Distributions of balanced accuracies from 100 Iris models over 4 hyperparameter sweeps.

A similar trend can be seen with MNIST models in Figure S.5. The balanced accuracies have tighter distributions than those seen in the entropy plots in Figure 4 of the main text, but they are not as apparent as those of the Iris models. We hypothesize that the differences we see in the location and size of the recurrent sets between the two models are also indicative of these observations. Namely, while both dynamical systems return two recurrent sets, the MNIST model's blue/green sets are near each other in Figure 5 while the shape of the recurrent sets is more complex for the Iris dataset in Figure 4.

Overall these trends indicate that models with different initial weights can achieve similar performance while not having a wide distribution of confidence in its predictions.

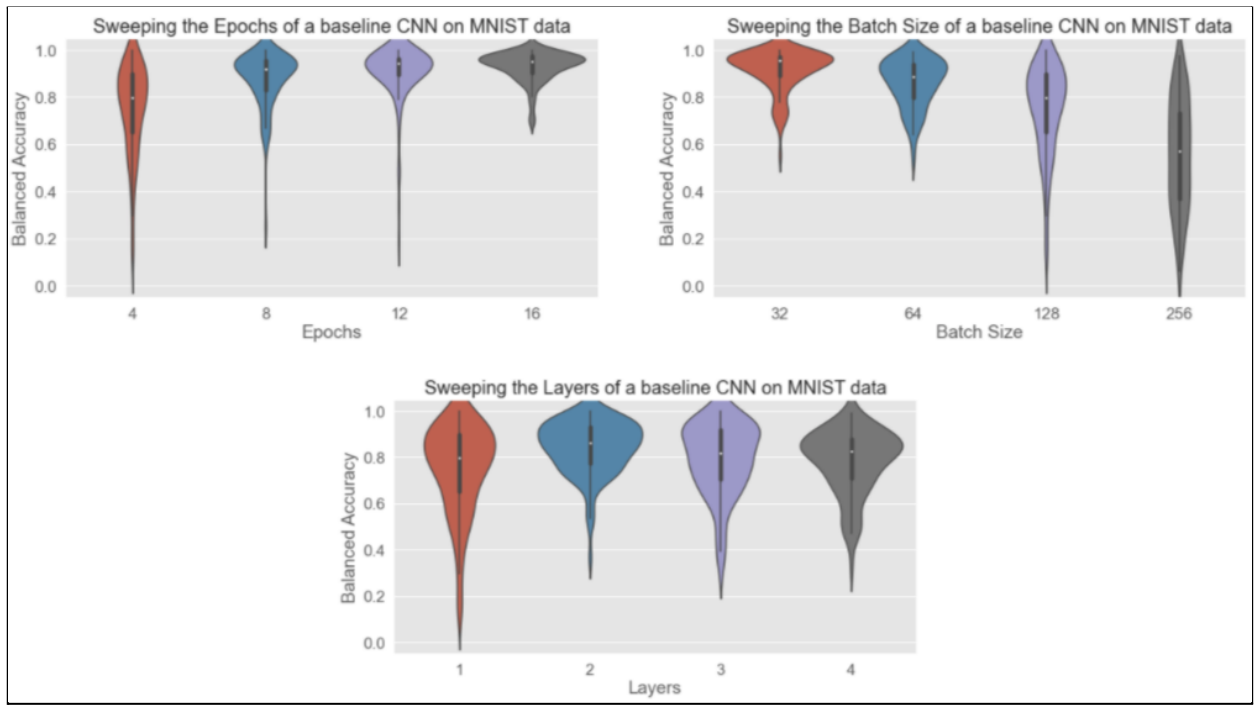

Figure S.5. Distributions of balanced accuracies from 100 MNIST models over 3 hyperparameter sweeps.

# S3. Selection of Weights to Summarize Model for Scaling CMGDB

CMGDB was run with different parameters of subdivisions and number of model parameters to analyze its scalability. As can be seen in Figure S.6, CMGDB run time increases exponentially as the number of phase subdivisions increases. Similarly, in Figure S.7 we see that CMGDB run time also increases exponentially as the number of weights increases. Thus, it was important for us to build a method to subselect the weights that exhibited the most "interesting" dynamics for further analysis.

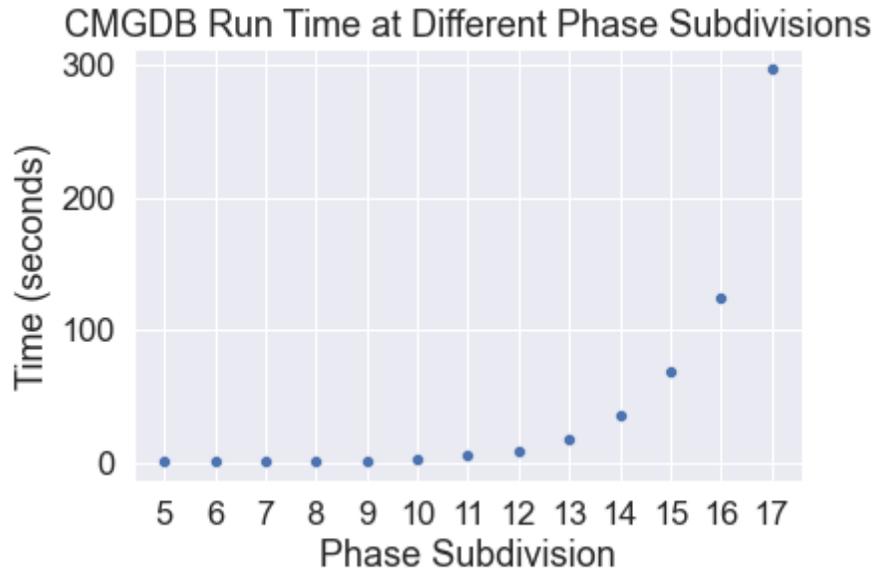

Figure S.6. CMGDB run time at different phase subdivisions. Data is from Iris Model with 450 Epochs.

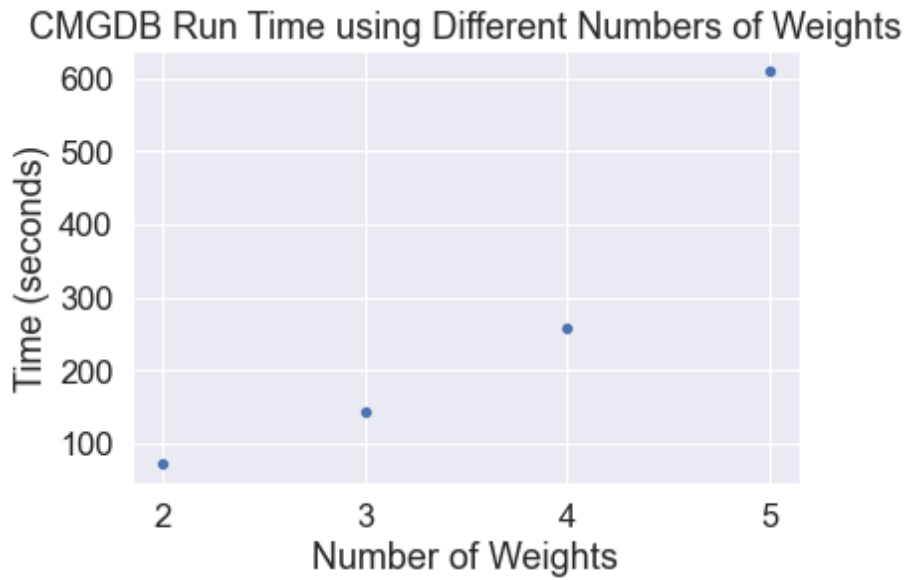

Figure S.7. CMGDB run time at different numbers of weights. Data is from Iris Model with 450 Epochs.

For each model architecture, the two weights with the highest variability between the trajectories were selected as "interesting." Namely, we excluded weights whose final state across the 100 instances were small (Figure S.8). The two weights of each model architecture were inspected for their variability over 100 trained models (Figure S.9) and the top two were selected for each instance.

6

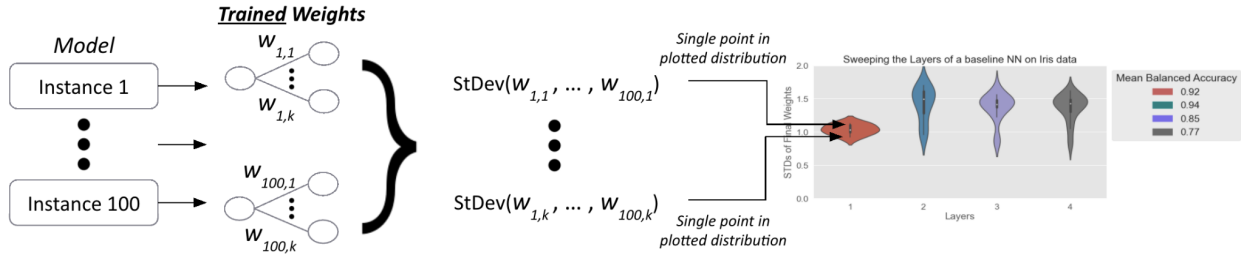

Figure S.8. illustration of how standard deviation values are calculated and used to create distributions of variability

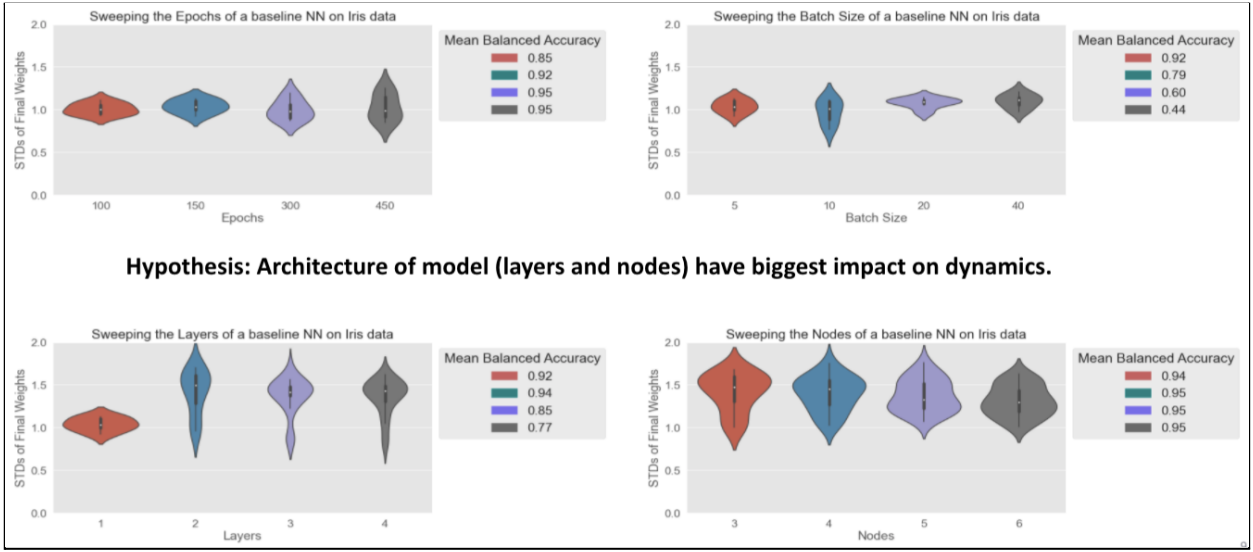

Figure S.9. Distributions of weight variability from 100 Iris models over 4 hyperparameter sweeps. Each point in each distribution represents the standard deviation of a weight over 100 trained models.

# S4. Numerical Construction of Multi-Valued Map, $\mathcal{F}$

Given a cubical grid, $\mathcal{X}$, we need to compute the image of grid elements by F and cover the image of a grid element by grid elements to determine the map $\mathcal{F}$. We compute $\mathcal{F}$ by applying the following procedure to each grid element (box) $\xi \in \mathcal{X}$:

1. Compute the image under F of the corner points of the grid element $\xi$.
2. Construct the smallest rectangle C (unrelated to the grid boxes) containing all the images of the corner points. Call this rectangle the rectangular cover of the image.
3. Find all grid boxes $\xi_1, \ldots, \xi_k$ in the cubical grid $\mathcal{X}$ that intersect the rectangular cover C of the image computed above.
4. Define $\mathcal{F}(\xi) = \{\xi_1, \ldots, \xi_k\}$, or as a digraph, add directed edges $\xi \to \xi_1, \ldots, \xi \to \xi_k$.

Applying the procedure above to all grid elements $\xi$ in $\mathcal{X}$ we obtain the multi-valued map $\mathcal{F}$.

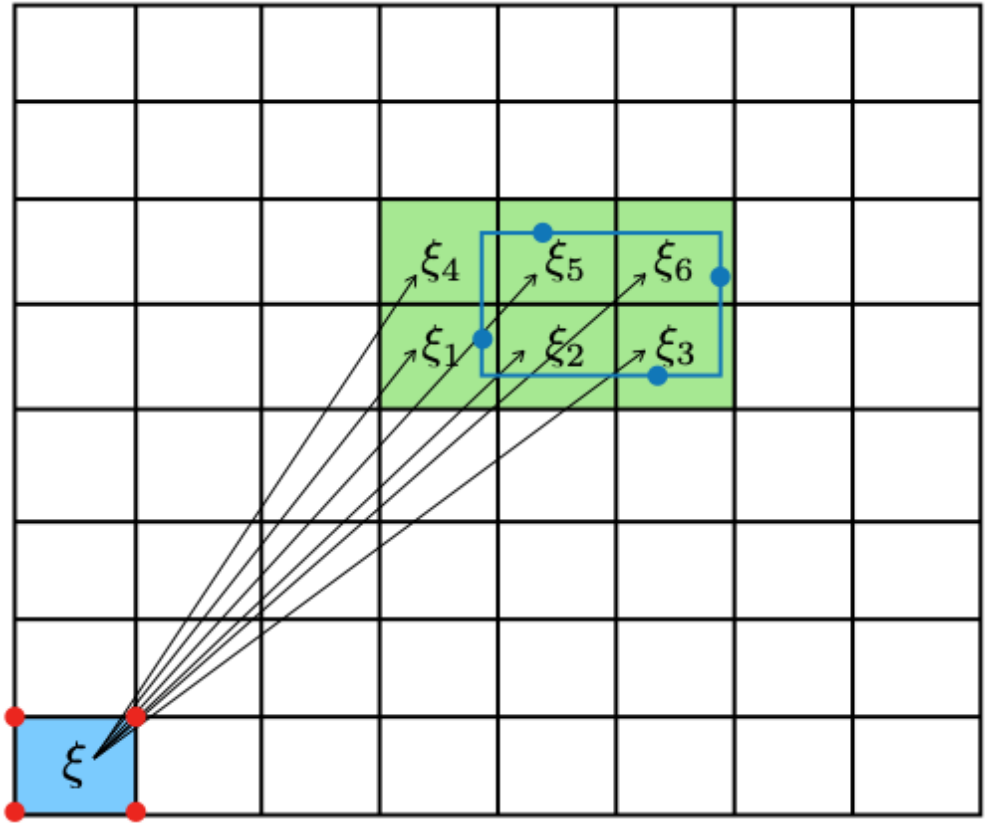

Figure S.10. Multi-valued map construction for a grid element ξ. The multi-valued map $\mathcal{F}$ maps the grid element ξ to the collection of grid elements {ξ$_1$, …, ξ$_6$}. The blue dots indicate the images of the corner points of ξ (red dots) under F. The green region indicates the elements intersecting the rectangle containing the blue dots.



\end{document}